\documentclass[12pt, reqno]{amsart}
\usepackage{amscd,amsthm,amsfonts,amssymb,amsmath,euscript}
\usepackage[dvips]{graphicx}
\usepackage[dvips]{graphics}
\usepackage[matrix,arrow]{xy}



\newcommand{\subvictor}{\par\medskip \noindent\refstepcounter{subsection}%
{\bf \arabic{section}.\arabic{subsection}.} }

\newcommand{\subsubvictor}{\par\medskip \noindent\refstepcounter{subsubsection}%
{\bf \arabic{section}.\arabic{subsection}.\arabic{subsubsection}.} }

\theoremstyle{definition}

\theoremstyle{remark}

\newcommand{\QQ}{{\mathbb Q}}
\newcommand{\ZZ}{{\mathbb Z}}
\newcommand{\NN}{{\mathbb N}}

\newcommand{\PP}{{\mathbb P}}
\newcommand{\CC}{{\mathbb C}}
\newcommand{\RR}{{\mathbb R}}
\newcommand{\BB}{{\mathbb B}}
\newcommand{\TT}{{\mathbb T}}

\newcommand{\Th}{{\bf Theorem}}

\newcommand{\Def}{{\bf Definition}}
\newcommand{\Rem}{{\bf Remark}}

\newcommand{\Corol}{{\bf Corollary}}

\newcommand{\Prop}{{\bf Proposition}}
\newcommand{\Prob}{{\bf Problem}}
\newcommand{\Proof}{{\bf Proof}}

\newfont{\smallskob}{cmbx7 scaled\magstep4}
\newfont{\bigskob}{cmbx12 scaled\magstep4}

\newcommand{\itc}[1]{\textup{#1}}

\newcommand{\pic}{\mathrm{Pic}\,}

\newcommand{\Div}{\mathrm{div}\,}
\newcommand{\Cl}{\mathrm{Cl}\,}

\newcommand{\Spec}{\mathrm{Spec}\,}
\newcommand{\Supp}{\mathrm{Supp}\,}
\newcommand{\Hom}{\mathrm{Hom}\,}

\newcommand{\gm}{{\mathbb G_{\mathrm{m}}}}

\newcommand{\tit}{On Landau--Ginzburg models for Fano varieties}

\begin{document}

\begin{title}
\tit
\end{title}

\begin{abstract}
We observe a method for finding weak Landau-Ginzburg models for Fano
varieties and find them for smooth Fano threefolds of genera $9$,
$10$, and $12$.
\end{abstract}

\author{Victor Przyjalkowski}

\thanks{The work was partially supported by FWF grant P20778,
    RFFI grants $08-01-00395$-a and $06-01-72017-$MNTI-a and grant
NSh$-1987.2008.1$}
\address{Steklov Mathematical Institute, 8 Gubkina street, Moscow 119991, Russia} %




\email{victorprz@mi.ras.ru, victorprz@gmail.com}

\maketitle

In the late 1980's physicists discovered a phenomenon of Mirror
Symmetry. They found that given a Calabi--Yau variety one may
construct the so called Superconformal Field Theory. This may be
done in two ways: ``algebro-geometric'' and ``symplectic''. Based on
this, they suggested that for each Calabi--Yau $X$ there is another
one $Y$ (not necessary uniquely determined), whose algebro-geometric
properties ``correspond'' to symplectic ones of $X$ and symplectic
ones ``corresponds'' to algebro-geometric ones of $X$. In
particular, the Hodge diamond of $Y$ is a reflection (a mirror
image) of the Hodge diamond of $X$ about a $45^\circ$ line. That's
why $X$ and $Y$ are called a mirror pair.

Later, in order to formalize this empiric approach, mathematicians
formulated a series of Mirror Symmetry Conjectures. They generalize
the correspondence to Fano varieties (Batyrev, Givental, Hori, Vafa,
etc.). The pair for a Fano variety $X$ is conjecturally a
Landau--Ginzburg model, that is, a (non-compact) manifold $M$ with
complex-valued function $f$ on it.

The dual model $M$ have a series of properties that correspond to
properties of $X$. Homological Mirror Symmetry Conjecture
(Kontsevich,~\cite{Ko94}), for instance, states that the derived
category of coherent sheaves on $X$ is isomorphic to the category of
Lagrangian vanishing cycles on $M$. This gives the correspondence
between singular fibers of $f$ and the exceptional collection of
$X$. One of the main problems in the Mirror Symmetry is to find a
Landau--Ginzburg model for a Fano variety. We use the Mirror
Symmetry Conjecture of Hodge structure variations to find a
candidate for it. This conjecture states the following. Given genus
0 two-pointed Gromov--Witten invariants of $X$ (the expected numbers
of rational curves of given degree that intersect two general
representatives of homological classes on $X$), one may define a
(small) Quantum Cohomology ring, i.~e. the deformation of the
cohomology of $X$ with Gromov--Witten invariants as structural
constants. The quantum multiplication in this ring gives the Quantum
$\mathcal D$-module. The conjecture is that the regularization of
this $\mathcal D$-module is isomorphic to a Picard--Fuchs $\mathcal
D$-module of $M$. The advantage of this conjecture is that we can
effectively check it in some cases. In the paper we check it in the
cases of some threefolds and discuss some approaches to finding
candidates to Landau--Ginzburg models. More about homological mirror
symmetry and mirror symmetry of Hodge structures variations for the
varieties we consider in the paper and for the other examples see
in~\cite{KP08}.

Consider a smooth Fano threefold $V$ with Picard group $\ZZ$. We may
associate a differential operator $L_V$ of type $D3$ with it
(see~\cite{Go05}, Definition $2.10$). The differential equation
associated with this operator has a unique analytical solution of
type $1+a_1t+a_2t^2+\ldots$ (the fundamental term of the regularized
$I$-series; see~\cite{Pr07} for the detailed description of the
solutions of equations of type $DN$). It is easy to calculate the
coefficients of the fundamental term using the recursive procedure.

Consider a Laurent polynomial $f\in
\CC[x,x^{-1},y,y^{-1},z,z^{-1}]$. Let $a_i$ be the constant term of
$f^i$.

\medskip

\Def. The polynomial $f$ is called \emph{a very weak
Landau--Ginzburg model} for $V$, if $a_i=b_i$ for all $i\geq 0$.

The polynomial $f$ is called \emph{a weak Landau--Ginzburg model}
for $V$, if $f$ is a very weak Landau--Ginzburg model for $V$ and a
general element of the pencil $\TT=\{1-tf=0, t\in \CC\}$ is
birational to a $K3$ surface.

\medskip

The similar definition may be formulated for the varieties of an
arbitrary dimension (see Definition~\ref{Definition: weak LG
models}).

The motivation of this definition is the following. The series
$1+b_1t+b_2t^2+\ldots$ is the solution of the Picard--Fuchs equation
for $\TT$ (see Proposition~\ref{Proposition: Picard--Fuchs}). The
Mirror Symmetry Conjecture of Hodge structures variations states
that the regularized Quantum Operator a variety (which is of type
$D3$ in the threefold case) coincides with the Picard--Fuchs
operator of its mirror dual Landau--Ginzburg model. The general
Mirror Symmetry philosophy says that the general fiber of the
Landau--Ginzburg model is a smooth Calabi--Yau variety.

\medskip

So, finding of weak Landau--Ginzburg models reduces to a
computational problem. That is, given $b_i$'s one should find $f$
with free terms $b_i$'s coinciding with $a_i$'s. The $K3$-condition
is usually follows from the degree reasons and the Bertini Theorem
(see below). The difficulty is computational. In the paper we
discuss some approaches to solving this problem. In particular, we
discuss conjectures related with toric degenerations of Fano
varieties. Using these approaches we find weak Landau--Ginzburg
models of Fano threefolds $V_{16}$ (of genus $9$), $V_{18}$ (of
genus $10$), and $V_{22}$ (of genus $12$).

\section{Definitions and conventions} \emph{The variety} is a
reduced irreducible scheme of finite type. The variety $X$ is called
\emph{$\QQ$-factorial} if $\Cl (X)\otimes \QQ\cong \pic (X)\otimes
\QQ$ (where $\Cl (X)$ is the group of Weil divisor classes on $X$).
It is said to be \emph{$\QQ$-Gorenstein} if $mK_X\in \pic(X)$ for
some $m\in \NN$. The $\QQ$-Gorenstein variety $X$ is said to have
\emph{canonical singularities} if for each resolution $f\colon
X^\prime \to X$ the relative canonical $\QQ$-divisor
$K_{X^\prime}-f^*(K_X)$ is effective.

\emph{The \itc{(}local\itc{)} deformation} is a flat morphism
$\mathcal X\to S$, where $S=(S,s_0)$ is a germ of a smooth variety
(usually the germ of a curve). The fiber over the central point
$X_{s_0}$ is called \emph{the special fiber}. The fiber over other
point is called \emph{the general fiber}. We say that the general
fiber \emph{degenerates} to the special fiber and $\mathcal X$ is
\emph{a degeneration} to $X_{s_0}$.

Let $X$ be a smooth algebraic variety with Picard group $\ZZ$. 
Let $\gamma_1,\ldots,\gamma_m\in H^*(X,\ZZ)$, $k_1,\ldots,k_m\in
\ZZ_{\ge 0}$, and $\beta\in H_2(X,\ZZ)$ be the class of algebraic
curves of anticanonical degree $d\geq 0$. The genus 0 Gromov--Witten
invariant with descendants that correspond to this data
(see~\cite{Ma99}, VI--2.1) we denote by
$$
\langle \tau_{k_1}
\gamma_1,\ldots,\tau_{k_m}\gamma_m\rangle_{\beta}=\langle \tau_{k_1}
\gamma_1,\ldots,\tau_{k_m}\gamma_m\rangle_{d}.
$$


Everything is over $\CC$.

\section{Weak Landau--Ginzburg models}

\subvictor Let $X$ be a smooth Fano variety of dimension $N$ with
Picard group $\ZZ$. Consider a torus
$\TT_{\mathrm{NS}^\vee}\cong\Spec B$, $B=\CC[t_0, t_0^{-1}]$, twice
dual to the Neron--Severi lattice. Let $H_H^*(X)\subset H^*(X,\QQ)$
be the subspace generated by the anticanonical class $H=-K_X$. This
subspace is tautologically closed with respect to the
multiplication, i. e. for any $\gamma_1, \gamma_2\in H^*_H(X)$ the
product $\gamma_1\cdot \gamma_2$ lies in $H_H^*(X)$. The
multiplication structure on the cohomology ring may be deformed.
That is, one can consider \emph{a quantum cohomology ring}
$QH^*(X)=H^*(X)\otimes \CC[t_0]$ (see~\cite{Ma99}, Definition
$0.0.2$) with quantum multiplication $\star\colon QH^*(X)\times
QH^*(X)\rightarrow QH^*(X)$, i.~e. the bilinear map given by
$$
\gamma_1\star \gamma_2=\sum_{\gamma,d} t_0^d\langle \gamma_1,
\gamma_2, \gamma^\vee\rangle_d \gamma,
$$
for all $\gamma_1,\gamma_2,\gamma\in H^*(X)$, where $\gamma^\vee$ is
the Poincar$\mathrm{\acute{e}}$ dual class to $\gamma$ (we identify
elements of $\gamma\in H^*(X)$ and $\gamma\otimes 1\in QH^*(X)$).
The constant term of $\gamma_1\star \gamma_2$ (with respect to
$t_0$) is $\gamma_1\cdot \gamma_2$. The subspace
$QH^*_H(X)=H^*_H(X)\otimes \CC[t_0]$ is not closed with respect to
$\star$ in general. The examples of varieties $V$ with non-closed
subspaces $H^*_H(X)$ are Grassmannians $G(k,n)$, $k, n-k>1$ of
dimension $>4$ (for instance, $G(2,5)$) or their hyperplane sections
of dimension $\ge 4$.

The variety is called \emph{quantum minimal} if $H^*_H(X)$ is closed
with respect to the quantum multiplication (see~\cite{Pr07},
Definition $1.2.1$). The examples of such varieties are Fano
complete intersections or Fano threefolds.

Let $HQ$ be a trivial vector bundle over $\TT_{\mathrm{NS}^\vee}$
with fiber $H_H^*(X)$. Let $S=H^0(HQ)$ and $\star\colon S\times S\to
S$ be the quantum multiplication (we may consider the quantum
multiplication as an operation on $S\cong QH^*_H(X)\otimes
\CC[t_0,t_0^{-1}]$).
Let $\mathcal D=B[t_0\frac{\partial}{\partial t_0}]$ and
$D=t_0\frac{\partial}{\partial t_0}$. Consider a \itc{(}flat\itc{)}
connection $\nabla$ on $HQ$ defined on the sections $H^i$ as
$$
\left( \nabla(H^i),t_0\frac{\partial}{\partial t_0}\right)=K_V\star
H^i
$$
(the pairing is the natural pairing between differential forms and
vector fields). This connection provides the structure of $\mathcal
D$-module for $S$ by $D(H^i)=(\nabla(H^i),D)$.

Let $Q$ be this $\mathcal D$-module. It is not regular in general.
To obtain the regular one, we need ``to regularize'' it. Let
$\gm=\Spec [t,t^{-1}]$. Let $E=\mathcal D_\gm/\mathcal D_\gm
(t\frac{\partial}{\partial t}-t)$ be the exponential $\mathcal
D_\gm$-module. Consider the inclusion $\ZZ(-K_X)\hookrightarrow \pic
X$. The natural isomorphism $\pic (X)\cong \mathrm{NS}(X)$ ($X$ is
Fano) and double dualization provide the morphism $j\colon \gm \to
\TT_{\mathrm{NS}^\vee}$. Define \emph{the regularization} of $Q$ as
$Q^{\mathrm{reg}}=j^*(\mu_*(Q\boxtimes j_*(E)))$, where $\mu\colon
\gm\times\gm\to \gm$ is the multiplication, and $\boxtimes$ is the
external tensor product (i.~e. $Q^{\mathrm{reg}}$ is a convolution
with the anticanonical exponential $\mathcal D$-module). It may be
represented as $\mathcal D_{\gm}/\mathcal
D_{\gm}(t\frac{\partial}{\partial t}L_X)$. We denote $t\frac{d}{dt}$
also by $D$. The differential operator $L_X$ is called \emph{of type
$DN$} (see~\cite{Go05}, 2.10). This operator is explicitly written
in~\cite{Go05}, Example $2.11$ for $N=3$ in terms of structural
constants of quantum multiplication by the anticanonical class (that
is,
two-pointed Gromov--Witten invariants). 
Thus, there is an operator of type $D3$ associated with every smooth
Fano threefold with Picard group $\ZZ$. There is $17$ families of
such Fanos (the Iskovskikh list). For all of them the operators of
type $D3$ are known, see, for instance,~\cite{Pr05}, $4.4$.

%

\subvictor \Def. { \itc{(}A unique\itc{)} analytic solution of
$L_XI=0$ of type
$$
I^X_{H^0}=1+a_1t+a_2t^2+\ldots\in \CC [[t]],\ \ \ \ b_i\in \CC,
$$
is called \emph{the fundamental term of the regularized $I$-series}
of $X$. }

\medskip

Let $\mathbf{1}$ be the class in $H^0(X,\ZZ)$ dual to the
fundamental class of $X$. Then this series is of the form
\begin{gather*}
I^X_{H^0}=1+\sum_{d\geq 2} \langle\tau_{d-2} \mathbf{1}\rangle_d
\cdot t^d
\end{gather*}
(see~\cite{Pr07}, Corollary $2.2.6$).

\subvictor Consider a torus $\TT=\mathbb
G_{\mathrm{m}}^n=\prod_{i=1}^n \Spec \CC[x_i,x_i^{-1}]$ and a
function $f$ on it. This function is represented by Laurent
polynomial: $f=f(x_1,x_1^{-1}\ldots,x_n,x_n^{-1})$. Let $\phi_f(i)$
be the constant term (i. e. the coefficient at $x_1^0\cdot \ldots
\cdot x_n^0$) of $f^i$. Put
$$
\Phi_f=\sum_{i=0}^\infty \phi_f(i)\cdot t^i\in \CC[[t]].
$$

\medskip

\Def. { The series $\Phi_f=\sum_{i=0}^\infty \phi_f(i)\cdot t^i$ is
called \emph{the constant terms series} of $f$. }

\subvictor \Def. { \label{Definition: weak LG models} Let $X$ be a
smooth $n$-dimensional quantum minimal Fano variety and
$I^X_{H^0}\in \CC[[t]]$ be its fundamental term of regularized
$I$-series. The Laurent polynomial $f\in \CC[\ZZ^n]$ is called
\emph{a very weak Landau--Ginzburg model} for $X$ if
$$
\Phi_f(t)=I^X_{H^0}(t).
$$

The Laurent polynomial $f\in \CC[\ZZ^n]$ is called \emph{a weak
Landau--Ginzburg model} for $X$ if it is a very weak
Landau--Ginzburg model for $X$ and for almost all $t\in \CC$ the
hypersurface $(1-tf=0)$ is birational to a Calabi--Yau variety. }

\medskip

The meaning of the definition is the following (see~\cite{SB85},
$10$, or~\cite{Be83}, pp. 50--52). Consider functions $F_t=1-t\cdot
f\in \CC[x_1,x_1^{-1},\ldots,x_n,x_n^{-1}][t]$. They provide a
pencil $\TT\to \BB=\PP[u:v]\setminus (0:1)$ with fibers
$Y_t=(F_t=0)$, $t\in \BB$. 

\medskip

The following proposition is a sort of a mathematical folklore.

\subvictor \Prop. {\it \label{Proposition: Picard--Fuchs} Let the
Newton polytope of $f\in \CC[\ZZ^n]$ contains $0$ in the interior.
Let $t\in \BB$ be the local coordinate around $(0:1)$. Then there is
a fiberwise $n-1$-form $\omega_t\in \Omega^{n-1}_{\TT/\BB}$ and
\itc{(}locally defined\itc{)} fiberwise $n-1$-cycle $\Delta_t$, such
that
$$
\Phi_f(t)=\int_{\Delta_t}\omega_t.
$$
}

\medskip

This means that $\Phi_f(t)$ is a solution of the Picard--Fuchs
equation for the pencil $\{Y_t\}$.

\medskip

\Proof. The following argument is based on~\cite{Gr69}, $\S 3$.

Let
$$
T_s=\{(x_1,\ldots,x_n)\in \TT: |x_1|=\ldots=|x_n|=s \},
$$
$T=T_1$, and $R_\delta=\cup T_s$, $\delta \leq s\leq 1$ (we consider
the natural metric on the torus given by $\TT\hookrightarrow
\CC^n$). Let $t$ be small enough such that $Y_t\cap T=\emptyset$,
that is, $|f(T)|<|1/t|$. Let $\delta$ be small enough, such that
$Y_t\cap T_\delta=\emptyset$ ($f$ has terms of negative degree). We
may assume that $R_\delta$ and $Y_t$ intersect transversally. Let
$\Delta_t=Y_t\cap R_\delta\in Y_t$. Let $Y^\varepsilon_t=\{p\in
\TT|\,\exists v\in Y_t: |p-v|<\varepsilon\}$, where $\varepsilon$ is
small enough. Then $\Delta^\varepsilon_t=R_\delta\cap
Y_t^\varepsilon$ is a ``tube'' over $\Delta_t$. Clearly
$\partial(R_\delta\setminus
Y^\varepsilon_t)=T+T_\delta-\Delta^\varepsilon_t$, so $T+T_\delta$
and $\Delta^\varepsilon_t$ are homologically equivalent.

Let
$$
\Omega_t=\frac{1}{(2\pi
i)^n}\frac{1}{F_t}\prod_{i=1}^n\frac{dx_i}{x_i}.
$$
Consider the integral
$$
\Phi (t)=\int_{T+T_\delta}\Omega_t.
$$
It is easy to see that $\int_{T_\delta}\Omega_t$ tends to zero under
$\delta\to 0$. Since it is constant, it equals zero. Therefore,
integrating step-by-step, we have $\Phi
(t)=\int_T\Omega_t=\Phi_f(t)$.

On the other hand, by Poincar$\rm \acute{e}$ Residue Theorem
$$
\Phi(t)=\int_{\Delta^\varepsilon_t}\Omega_t=\int_{\Delta_t}{\rm
Res}_{Y_t} \Omega_t=\int_{\Delta_t}\omega_t.
$$
\qed

\subvictor \label{remark:equations} Let $PF_f=PF_f(t,
\frac{\partial}{\partial t})$ be a Picard--Fuchs operator of
$\{Y_t\}$. Let $m$ be the order of $PF_f$ and $r$ be the degree with
respect to $t$. Let $Y$ be a semistable compactification of
$\{Y_t\}$ (so we have the map $\widetilde{f}\colon Y\to \PP^1$;
denote it for the simplicity by $f$). Let $m_f$ be the dimension of
transcendental part of $R^{n-1}f_!\,\ZZ_Y$ (the algorithm for
computing it see in~\cite{DH86}). Let $r_f$ be the number of
singularities of $f$ (counted with multiplicities). Then, $m\leq
m_f$ and $r\leq r_f$. Thus, the first few coefficients of the
expansion of the solution of the Picard--Fuchs equation determine
the other ones. This means that if the first few coefficients of the
expansion the solution of $L\Phi=0$, $L\in \CC[t,
\frac{\partial}{\partial t}]$, coincide with the first few
coefficients of the expansion of $\Phi_f$, then $L=PF_f$.


\section{Main theorem}

\subvictor \Th. {\it \label{Theorem: main}
\begin{enumerate}
    \item The Laurent polynomial
\begin{multline*}
f_{16}=\frac{1}{xyz}+2\left(\frac{1}{xy}+\frac{1}{xz}+
\frac{1}{yz}\right)+3\left(\frac{1}{x}+\frac{1}{y}+\frac{1}{z}\right)+\left(\frac{x}{yz}+\frac{y}{xz}+\frac{z}{xy}\right)\\
\left(\frac{x}{y}+\frac{x}{z}+\frac{y}{x}+\frac{y}{z}+\frac{z}{x}+\frac{z}{y}\right)+4+\left(x+y+z\right)
\end{multline*}
is a weak Landau--Ginzburg model for Fano variety $V_{16}$.
    \item The Laurent polynomial
\begin{multline*}
f_{18}=2\left(\frac{1}{x}+\frac{1}{y}+\frac{1}{z}\right)+\left(\frac{x}{yz}+\frac{y}{xz}+\frac{z}{xy}\right)+\\
\left(\frac{x}{y}+\frac{x}{z}+\frac{y}{x}+\frac{y}{z}+\frac{z}{x}+\frac{z}{y}\right)+3+\left(x+y+z\right)
\end{multline*}
is a weak Landau--Ginzburg model for Fano variety $V_{18}$.
    \item The Laurent polynomial
$$
f_{22}=
\frac{xy}{z}+\frac{y}{z}+\frac{x}{z}+x+y+\frac{1}{z}+4+\frac{1}{x}+\frac{1}{y}+z+\frac{1}{xy}+\frac{z}{x}+\frac{z}{y}+
\frac{z}{xy}
$$
is a weak Landau--Ginzburg model for Fano variety $V_{22}$.
\end{enumerate}
}

\medskip

\Proof. The operators of type $D3$ are
$$
D^3-4t(2D+1)(3D^2+3D+1)+16t^2(D+1)^3
$$
for $V_{16}$,
$$
D^3-3t(2D+1)(3D^2+3D+1)-27t^2(D+1)^3
$$
for $V_{18}$, and
\begin{multline*}
D^3-\frac{2}{5}t(2D+1)(17D^2+17D+16)-\frac{56}{25}t^2(D+1)(11D^2+22D+12)-\\
\frac{126}{125}t^3(D+1)(D+2)(2D+3)-\frac{1504}{625}t^4(D+1)(D+2)(D+3)
\end{multline*}
for $V_{22}$ (see~\cite{Go05}). The degrees of Picard--Fuchs
operators for pencils that are given by $f_{16}$, $f_{18}$, and
$f_{22}$ are bounded by $3$ with respect to $D$ and $4$ with respect
to $t$. One may check that the first few coefficients of the
expansion of the constant terms series and the fundamental term of
regularized $I$-series coincide (see~\ref{remark:equations}). Thus,
these polynomials are very weak Landau--Ginzburg models.

Compactify the torus to $\PP^3$ in the standard way. Then the
elements of the pencils that are given by $f_{16}$, $f_{18}$, and
$f_{22}$ are quartics in $\PP^3$, so, have trivial canonical class.
By Bertini's Theorem, the general element of the pencil may have
only base points as singularities. It is easy tho check, that all
base points of any pencil are Du Val, so the canonical class of
minimal model of the general fiber of it is trivial. By~\cite{Ch96},
Theorem $5.1$, the general fiber is birational to a $K3$ surface.
\qed

\subvictor \Rem. The definition of weak Landau--Ginzburg model is
numerical. So, it is natural that the polynomials from
Theorem~\ref{Theorem: main} are not unique. For instance,
polynomials obtained from them by some coordinate changes or
resizing $x\to \alpha x$, $y\to \beta y$, $z\to \gamma z$ are also
weak Landau--Ginzburg models.

%

%

\subvictor {\bf Singularities.} The philosophy of Mirror Symmetry
says that for every smooth variety there is a dual Landau--Ginzburg
model, that is, the pencil of projective varieties, whose symplectic
properties correspond to algebro-geometric ones of the variety and
vice-versa. In particular, the sheaf of vanishing cycles on the
element of the pencil ``lying over 0'' corresponds to the horizontal
Hodge cohomologies of the variety and the elements with isolated
singularities correspond to the bounded derived category of coherent
sheaves on the variety. The fibers of weak Landau--Ginzburg
models are non-compact. 
%
Compactify these models (as lying in $\TT\hookrightarrow\PP^3$) and
resolve singularities of the compactifications. Then the
singularities of the fibers are\footnote{Remind that the fiber of
weak Landau--Ginzburg model $f$ over infinity, which is the fiber
over 0 in the standard coordinates, is given by $f=0$.}:
\begin{description}
    \item[$f_{16}$] the (reducible) curve of genus three over the infinity 
    and two conjugate points defined over the quadratic extension of $\QQ$.
    \item[$f_{18}$] the (reducible) curve of genus two over the infinity and two
    conjugate points defined over the quadratic extension of $\QQ$.
    \item[$f_{22}$] 
    three conjugate points defined over the extension of $\QQ$ of degree 3.
\end{description}
The fibers of the pencils over the infinity in these cases coincide
with the expectation. Namely, the genera of the schemes of
singularities equal the dimensions of the intermediate Jacobians of
the corresponding varieties.

The images of singular points are the singular points of the
differential operators of type $D3$ for these varieties. More about
these models see in~\cite{KP08}.

\subvictor \Rem. The aim of this paper is to describe
Landau--Ginzburg models for smooth Fano threefolds with Picard group
$\ZZ$. For most of them (for complete intersections in projective or
weighted projective varieties) they are known (see~\cite{HV00},
$7.2$). The last cases that we have not found yet
are $V_{10}$, $V_{12}$, and $V_{14}$.

%
%

\section{Finding weak Landau--Ginzburg models}

Unfortunately, finding weak Landau--Ginzburg models is very
complicated computational problem. In this section we discuss some
(empiric) ways of simplifying it.

\subvictor {\bf Canonical degenerations and numerical invariants}

\subsubvictor \Th\ [Kawamata,~\cite{Ka97}]. {\it Let $\mathcal X\to
S$ be a deformation \itc{(}$S$ is a germ of a curve\itc{)}. Suppose
that $X_{s_0}$ has canonical singularities. Then $\mathcal X$ has
canonical singularities. In particular, for any $s\in S$ the fiber
$X_s$ is canonical. }

\subsubvictor \Corol. {\it \label{Corollary: degree} The total space
$\mathcal X$ is $\QQ$-Gorenstein. Thus, by adjunction, for any $s\in
S$ we have $-K_{X_s}=-K_{\mathcal X}|_{X_s}$. In particular, the
anticanonical degree $(-K_s)^{\dim X_s}$ does not depend on $s\in
S$. }

\subsubvictor \label{Statement: genera} Let $\mathcal F$ be a
coherent sheaf on $\mathcal X$ which is flat over $S$. Then the
Euler characteristic $\chi (X_s,\mathcal F_s)$ does not depend on
$s\in S$ (see, for instance,~\cite{Da96}, Proposition 3.8). Let
$X_s$ be a canonical Fano variety for every $s\in S$, $s\neq s_0$
and $X_{s_0}$ be a canonical almost Fano variety (that is, its
anticanonical divisor is
nef and big). 
By Kawamata--Viehweg Vanishing Theorem (Theorem 2.17
in~\cite{Ko96}), $\chi (X_s,-K_s)=h^0(X_s,-K_s)$. Thus,
$h^0(X_s,-K_s)$ does not depend on $s\in S$.

%

\subsubvictor \Prop. {\it \label{Proposition: Picard rank} Let
$\pi\colon\mathcal X\to S$ be a deformation such that the general
fiber is a Fano variety of Picard rank $k$ and the special fiber
$X_{s_0}$ is irreducible, projective, and normal almost Fano
variety. Let all fibers have canonical singularities.
Let $\pic (X_{s_0})=\ZZ^m$. Then $m\leq k$. }

\medskip

\Proof\ (the idea is due to Ivo Radloff). We use here the Picard
groups and cohomology groups with coefficients in $\QQ$. Let
$\Delta=\{t:|t-s_0|<\varepsilon\}\subset S$ be a small enough
neighbourhood of $s_0$ and $X=\pi^{-1} (\Delta)$. Then $H^2 (X)=H^2
(X_{s_0})$ as $X_{s_0}$ is a deformation retract of $X$. By
Kawamata--Viehweg Vanishing Theorem and an exponential exact
sequence
$H^2(X)\cong \pic (X)$ and $H^2 (X_{s_0})\cong \pic (X_{s_0})$ (
$X$ is a relative Fano variety). Thus, we need to show that there is
no linear sheaf $\mathcal L$ such that the restriction $\mathcal
L|_{X_s}\cong \mathcal O_{X_s}$ for $s\neq s_0$ and positive for
$s=s_0$ (the numerical equivalence over $\QQ$ is the same as the
linear equivalence).

Suppose it is. By semicontinuity (see, for example, Theorem $3.6$
in~\cite{Da96}), there is a section of $\mathcal L|_{s_0}$. It is
non-zero by assumption, so it is an effective divisor. Denote the
dimension of the fibers by $n$. (We apply an intersection theory to
sheaves as to the linear systems associated to them.) The special
fiber is projective, so there is a divisor $\mathcal D$ on $X$ whose
restriction on the special fiber is ample. So, $\mathcal
D^{n-1}\cdot \mathcal L\cdot X_{s_0}=(\mathcal
D|_{X_{s_0}})^{n-1}\cdot \mathcal L|_{s_0}>0$. The intersection
number does not depend on the fiber as all fibers are numerically
equivalent. The sheave $\mathcal L$ restricted to the general fiber
is numerically trivial by assumption. Contradiction. \qed

\subsubvictor \Corol. {\it \label{Corollary: invariants}
Let $\mathcal X\to S$ be a degeneration of 
Fano variety of Picard rank $1$ to the toric canonical variety
$X_{s_0}$. Then the anticanonical degree $(-K_{X_s})^{\dim X_s}$,
$h^0(-K_{X_s})$, and the Picard rank of $X_s$ do not depend on $s\in
S$. }

\medskip

 \Proof.
It follows from \ref{Corollary: degree}, \ref{Statement: genera},
and \ref{Proposition: Picard rank}. \qed

\subvictor {\bf Toric varieties and Laurent polynomials.} Consider a
torus $\TT=\Spec \CC[M]\cong \Spec (\CC^*)^n$, where $M\cong \ZZ^n$.
Let $x_1,\ldots,x_n$ be the coordinates on the one-dimensional tori.
Put $x^m=x_1^{m_1}\cdot\ldots \cdot x_n^{m_n}$ for
$m=(m_1,\ldots,m_n)$. Then any function $f$ on $\TT$ can be uniquely
represented as $f=\sum_{m\in M} a_m x^m$. Let $\Supp(f)=\{m\in M,
a_m\neq 0\}$.
The convex hull of $\Supp(f)$ in $M_{\RR}=M\otimes \RR$ is called
\emph{the Newton polyhedra} of $f$.
%

Let $X$ be a toric variety with an open subset $\TT$. Let $N=\Hom
(M,\ZZ)$ be the dual to $M$ lattice and $\langle \cdot,\cdot\rangle$
be a natural pairing. It is well-known (see, for
instance,~\cite{Da78}) that $X$ is associated with \emph{the fan}
$\Sigma\subset N_{\RR}=N\otimes \RR$. The variety $X=X_\Sigma$ is
covered by affine toric varieties (maps)
$X_\sigma=\Spec\CC[\sigma^\vee\cap M]$, where $\sigma^\vee=\{m\in
M_\RR|\, \langle m,v\rangle\geq 0, v\in\sigma\}$ are dual cones for
$1$-dimensional cones $\sigma\in\Sigma$.

In the following we assume that all fans are \emph{projective}, that
is, they correspond to the projective varieties (this is equivalent
to the existence of strictly convex function on the fan which is
linear on all cones from $\Sigma$). All polytopes are supposed to be
convex and to have the origin in the interior.

Let $L\cong \ZZ^n$ be any lattice, $L_\QQ=L\otimes \QQ$,
$L_\RR=L\otimes \RR$, and $L^\vee$ be its dual. For any fan
$\Sigma\in L$ we associate the polytope $P_\Sigma\in L_\RR$ which is
defined to be a convex hull of primitive vectors of its rays (that
is, primitive vectors $v_i\in L$ that generate one-dimensional cones
of
$\Sigma$). 
Otherwise,
given any polytope $P\in L_\QQ$, we can construct a \emph{normal
fan} 
taking cones over its faces. Given a polytope in $L_\RR$, define its
dual as
$$
P^\vee=\{ m\in L^\vee_\RR| \langle m,n\rangle\geq -1 \mbox{ for all
} n\in P \}.
$$
Obviously, if $P$ have vertices in $L_\QQ$, then $P^\vee$ have
vertices in $L^\vee_\QQ$. Therefore, given $P$ we may construct a
toric variety $X_P$ that is given by a normal fan for
$P^\vee$\footnote{In fact, the datum $(P,\varphi)$, where $P\in L$
is an integral polytope containing origin and $\varphi$ is a
strictly convex integral piecewise linear function is equivalent to
the fan in $L^\vee$. This function here is the function $\varphi$
such that $\varphi(n)=-1$ for all vertices of $P$ (it is represented
by the anticanonical class).}. The polytope $P\subset L_\RR$ with
vertices in $L$ is called \emph{reflexive} if $P^\vee$ have vertices
in $L^\vee$. The toric variety that correspond to a reflexive
polytope is a Gorenstein Fano with at most canonical singularities.

Let $\Sigma\subset N$ be a fan, $P=P_\Sigma\subset N_\RR$ be the
corresponding polytope and $P^\vee\subset M_\QQ$ be its dual. The
anticanonical divisor of $X_\Sigma$ is the sum of boundary divisors
$D_1,\ldots,D_r$ corresponding to rays given by primitive vectors
$n_1,\ldots,n_r$. The point $m\in M$ is a rational function on
$X_\Sigma$. Its divisor is $\sum \langle m, n_i\rangle D_i$. Let us
consider $\QQ$-divisors in the following. The element of $M_\QQ$
determines a $\QQ$-divisor by linearity. In particular, the Newton
polyhedra $\Delta\in M$ of the function $f$ lies in $P^\vee$ if and
only if $\Div (f)-K_{X_\Sigma}\in \pic(X_\Sigma)\otimes \QQ$ is
effective (where $\Div (f)$ is the divisor of $f$). Thus, functions
whose Newton polytope lie in $P^\vee$ are the sections of the
anticanonical sheaf, so $L(P^\vee)$
is naturally isomorphic to $|-K_{X_\Sigma}|$, where $L(P^\vee)$ is
the space of Laurent polynomials with support in $P^\vee$.

Let $P\subset N_\RR\cong \ZZ^n\otimes \RR$ be a polytope and $X$ be
the toric variety associated with the normal fan of $P$. Then $X$
has canonical singularities if and only if $P$ has the origin as the
only integral point in the interior. The the anticanonical degree of
$X$ (that is, $(-K_X)^n$) is a volume of the dual polytope $P^\vee$
divided by $n!$. The dimension of the anticanonical linear system
equals the number of integral points inside $P^\vee$ and on the
boundary. The Picard number of $X$ equals the dimension of the space
of functions on $N_\RR$ that are linear on every cone over the face
of $P$ modulo linear functions.

\subvictor {\bf Strategy.}




So, the natural strategy for finding weak Landau--Ginzburg models of
Fano threefolds is the following. Consider a smooth Fano threefold
$X$ with Picard group $\ZZ$. Find the fundamental term
$I^X_{H^0}=\sum a_rt^r$, $b_r\in \QQ$, of its $I$-series (these
series are known for all 17 families of such Fanos, see, for
instance,~\cite{Pr05}, 4.4). Suppose that $X$ degenerates to a
canonical Fano $X_\Sigma$ and suppose that there exist a weak
Landau--Ginzburg model $f$ for $X$, whose Newton polyhedra lies in
$P_{\Sigma}$. Find it. For this find all integral polytopes with the
origin as a unique integral point in the interior, whose numerical
data is the same as the data of $X$. Consider any such polytope and
the Laurent polynomial $f=\sum b_{ijk}x^iy^jz^k\in
\CC[x,x^{-1},y,y^{-1},z,z^{-1}][b_{ijk}]$, whose Newton polyhedra is
our polytope. Let $\Phi_f=\sum b_r(b_{ijk}) t^r$ be its constant
terms series. To specify the coefficients $b_{ijk}$, solve the
system of equations $\{b_r(b_{ijk})=a_r\}$, $r=1,\ldots, N$, where
$N\in \NN$ is big enough; to avoid resizing, normalize $x,y$ and $z$
such that $b_{100},b_{010}$ and $b_{001}$ are $0$ or $1$.
Prove that $\Phi_f=I^X_{H^0}$ for $f$ we found. For this check that
for all coefficients of $\Phi_f$ holds the same recurrence that
holds for the coefficients of $I^X_{H^0}$. Finally, to prove that
the general element of our pencil is birational to a Calabi--Yau
variety it usually suffices to compactify the torus to the
projective space, compactify the fibers therein, check the degree
condition, use the Bertini Theorem and check that the general
hypersurface admits a crepant resolution.

\medskip

To ``legalize'' this empiric strategy, one should solve two
following problems.

\subsubvictor \Prob. {\it Prove that any smooth Fano threefold with
Picard group $\ZZ$ admits a degeneration to a canonical toric
variety. Find all such degenerations. Characterize them. Generalize
this to more general class of Fanos or to toric varieties with worse
singularities.}

\medskip

This problem may be solved if the singularities of toric variety are
terminal Gorenstein (that is, ordinary double points)
(see~\cite{Ga07}). Unfortunately, there is only five such varieties,
that is, $\PP^3$, three-dimensional quadric, complete intersection
of two quadrics, and the varieties $V_5$ and $V_{22}$. It is
remarkable that we do not need the particular form of degeneration
(but in some cases, as for quadric or complete intersection of two
quadrics we can find them, see~\cite{Ba97}). Unfortunately, we can't
expect that any smooth Fano variety degenerates to a
\emph{Gorenstein} canonical toric Fano variety (that is, associated
with a reflexive polytope). The example is $V_2$, the Fano variety
of degree $2$: there is no reflexive polytopes of volume
$\frac{1}{3}$.

\subsubvictor \Prob. {\it Let the smooth Fano variety $X$ degenerate
to the canonical toric variety $T$. Prove that there is a weak
Landau--Ginzburg model $f$ for $X$ with Newton polytope $\Delta$ and
there is a fan $\Sigma$ of $T$ such that $P_\Sigma=\Delta^\vee$. }

\medskip

The good references for this problem are~\cite{Ba97}
and~\cite{BCFKS98a}. In the paper~\cite{EHX97} the first examples of
weak Landau--Ginzburg mirrors for nontoric Fano varieties with
Picard number 1 (Grassmannians) were suggested. In the
paper~\cite{BCFKS98b} their relation to toric degenerations of
Grassmannians was explained.

\subsubvictor Unfortunately, this straightforward way is too
complicated for computational reasons. First, there are too many
such polytopes. Secondly, there are too many integral points in any
three-dimensional polytope, so there are too many variables in the
system of equations. Thirdly, it is complicated to solve the system
of polynomial equations.

To fix these problems we put some restrictions.

To fix the first problem we consider not all such polytopes, but
some natural class of them (such as reflexive polytopes, polytopes
with many symmetries or polytopes that are contained in the cube
$[-1,1]\times [-1,1]\times [-1,1]$) and hope that $X$ degenerates to
the toric variety associated with such polytope. As the degree of
Fano variety increases, the polytope of its degeneration tends to
become simpler.

To fix the second one we consider not all functions but functions of
some type. Namely, we may consider only Laurent polynomials
$f=f(x,y,z)$ which is symmetric under all permutations of $x,y,z$.
We may also consider polynomials with coefficients $1$ at the
vertices of their Newton polytopes. Polynomials we found are of
these types.

Finally, we hope that the coefficients of polynomials are integral.
So, to solve the system of equations we solve it modulo some prime
numbers, lift the solutions to $\ZZ$ and check if we did this
correct. Actually, we consider all possibilities for
$b_{ijk}\!\!\mod p$ and check if the equations hold for them.

\bigskip

The author is grateful to V.\,Golyshev for proposing the problem and
explanations, I.\,Radloff for the idea of the proof of
Proposition~\ref{Proposition: Picard rank}, L.\,Katzarkov,
Yu.\,Prokhorov, K.\,Shramov, and the referee for references on
papers~\cite{EHX97} and~\cite{BCFKS98b}.

\end{document}